\documentclass[11pt]{article}
\usepackage{eurosym}
\usepackage{color}
\usepackage{amsmath}
\usepackage{amsfonts}
\usepackage{amssymb}
\usepackage{authblk}
\usepackage{amsthm}
\usepackage{dsfont}
\usepackage{mathtools}

\theoremstyle{definition}

\begin{document}

\title{Tanaka's Theorem Revisited}
\author{Saeideh Bahrami \\
{\small {bahrami.saeideh@gmail.com}}}
\maketitle

\begin{abstract}
Tanaka (1997) proved a powerful generalization of Friedman's self-embedding
theorem that states that given a countable nonstandard model $(\mathcal{M},%
\mathcal{A})$ of the subsystem $\mathrm{WKL}_{0}$ of second order
arithmetic, and any element $m$ of $\mathcal{M}$, there is a self-embedding $%
j$ of $(\mathcal{M},\mathcal{A})$ onto a proper initial segment of itself
such that $j$ fixes every predecessor of $m$.

Here we extend Tanaka's work by establishing the following results for a
countable nonstandard model $(\mathcal{M},\mathcal{A})\ $of $\mathrm{WKL}%
_{0} $ and a proper cut $\mathrm{I}$\ of $\mathcal{M}$: \medskip

\noindent \textbf{Theorem A.}~\textit{The following conditions are equivalent%
}:

\noindent \textbf{(a)} $\mathrm{I}$\ \textit{is closed under exponentiation}.

\noindent \textbf{(b) }\textit{There is a self-embedding} $j$ \textit{of} $(%
\mathcal{M},\mathcal{A})$ \textit{onto a proper initial segment of} \textit{%
itself} \textit{such that }$I$\textit{\ is the longest initial segment of
fixed points of} $j$.\medskip

\noindent \textbf{Theorem B.}~\textit{The following conditions are equivalent%
}:

\noindent \textbf{(a)} $\mathrm{I}$\ \textit{is a strong cut of }$\mathcal{M}
$\textit{\ and }$\mathrm{I}\prec _{\Sigma _{1}}\mathcal{M}.$

\noindent \textbf{(b) }\textit{There is a self-embedding} $j$ of $(\mathcal{M%
},\mathcal{A})$ \textit{onto a proper initial segment of itself such that }$%
\mathrm{I} $\textit{\ is the set of all fixed points of} $j$\textit{.}
\end{abstract}

\section{Introduction}

One of the fundamental results concerning nonstandard models of Peano
arithmetic $(\mathrm{PA}\mathsf{)}$ is Friedman's theorem \cite[Theorem 4.4]%
{fri se} that states every countable nonstandard model of $\mathrm{PA}$ is
isomorphic to a proper initial segment of itself. A notable generalization
of Friedman's theorem was established by Tanaka \cite{tan} for models of the
well-known subsystem $\mathrm{WKL}_{0}$ of second order arithmetic, who
established: \medskip

\noindent \textbf{Theorem 1.1.}~\textbf{(Tanaka)} \textit{Suppose }$(\mathcal{M},%
\mathcal{A})$ \textit{is a countable nonstandard model} \textit{of }$\mathrm{%
WKL}_{0}$\textit{. }

\noindent \textbf{(a)} \textit{There is a proper initial segment }$\mathrm{I}
$\textit{\ of }$\mathcal{M}$\textit{\ and an isomorphism }$j$\textit{\
between }$(\mathcal{M},\mathcal{A)}$ \textit{and} $(\mathrm{I},\mathcal{A}_{%
\mathrm{I}})$\textit{, where }$\mathcal{A}_{\mathrm{I}}:=\{\mathrm{A}\cap
\mathrm{I}:\ \mathrm{A}\in \mathcal{A}\}$\textit{.}

\noindent \textbf{(b)} \textit{Given any prescribed }$m$ \textit{in} $%
\mathcal{M}$\textit{, there is an }$\mathrm{I}$\textit{\ and }$j$\textit{\
as in (a) such that }$j(x)=x$\textit{\ for all }$x\leq m.$\textit{\medskip }

Tanaka's principal motivation in establishing Theorem 1.1 was the
development of non-standard methods within $\mathrm{WKL}_{0}$ in the context
of the reverse mathematics research program; for example Tanaka and Yamazaki
\cite{haar} used Theorem 1.1 to show that the Haar measure over compact
groups can be implemented in $\mathrm{WKL}_{0}$ via a detour through
nonstandard models. This is in contrast to the previously known
constructions of the Haar measure whose implementation required the stronger
subsystem \textrm{ACA}$_{0}$. Other notable applications of the methodology
of nonstandard models can be found in the work of Sakamato and Yokoyama \cite%
{jor}, who showed that over the subsystem $\mathrm{RCA}_{0}$ the Jordan
curve theorem and the Sch\"{o}nflies theorem are equivalent to $\mathrm{WKL}%
_{0}$; and in the work of Yokoyama and Horihata \cite{rim}, who established
the equivalence of $\mathrm{ACA}_{0}$ and Riemann's mapping theorem for
Jordan regions over $\mathrm{WKL}_{0}.$\textit{\medskip }

Here we continue our work \cite{our} on the study of fixed point sets of
self-embeddings  of countable nonstandard models of I$\Sigma _{1}$
by focusing on the behavior of fixed point sets in Tanaka's theorem. Our
methodology can be generally described as an amalgamation of Enayat's
strategy for proving Tanaka's theorem \cite{Me Tanaka} with some ideas and
results from \cite{our}.\footnote{%
Enayat's paper \cite{Me Tanaka} provides a complete proof of part (a) of
Theorem 1.1, and an outline of the proof of part (b)\ of Theorem 1.1.} \
Before stating our results, recall that $j$ is said to be a proper initial
self-embedding of $(\mathcal{M},\mathcal{A)}$ if there is a proper initial
segment $\mathrm{I}$\ of\textit{\ }$\mathcal{M}$\textit{\ }and an isomorphism%
\textit{\ }$j$\textit{\ }between\textit{\ }$(\mathcal{M},\mathcal{A)}$ and $(%
\mathrm{I},\mathcal{A}_{\mathrm{I}})$\textit{, }where\textit{\ }$\mathcal{A}%
_{\mathrm{I}}:=\{\mathrm{A}\cap \mathrm{I}:\ \mathrm{A}\in \mathcal{A}\}$%
\textit{; }$\mathrm{I}_{\mathrm{fix}}(j)$ is the longest initial segment of
fixed points of $j$; and $\mathrm{Fix}(j)$ is the fixed point set of $j$, in
other words:

\begin{center}
$\mathrm{I}_{\mathrm{fix}}(j):=\{m\in M:\forall x\leq m\ j(x)=x\},$ and\\[0pt%
]
$\mathrm{Fix}(j):=\{m\in M:j(m)=m\}.$
\end{center}

Our main results are Theorems A and B below. Note that Theorem A is a
strengthening of Tanaka's Theorem (see Section 3 for more detail). \textit{%
\medskip }

\noindent \textbf{Theorem A.}~\textit{Suppose} $(\mathcal{M},\mathcal{A})$
\textit{is a countable nonstandard model of} $\mathrm{WKL}_{0}.$\textit{The
following conditions are equivalent for a proper cut} $\mathrm{I}$ \textit{of%
} $\mathcal{M}$:\medskip

\noindent \textbf{(1)} \textit{There is a self-embedding} $j$ \textit{of} $(%
\mathcal{M},\mathcal{A})$ \textit{such that }$\mathrm{I}_{\mathrm{fix}}(j)=%
\mathrm{I}. $\medskip

\noindent \textbf{(2) }$\mathrm{I}$\ \textit{is closed under exponentiation}%
.\medskip

\noindent \textbf{(3) }\textit{There is a proper initial self-embedding} $j$
\textit{of} $(\mathcal{M},\mathcal{A})$ \textit{such that }$\mathrm{I}_{%
\mathrm{fix}}(j)=\mathrm{I}.$\textit{\medskip }

\noindent \textbf{Theorem B.}~\textit{Suppose} $(\mathcal{M},\mathcal{A})$
\textit{is a countable nonstandard model of} $\mathrm{WKL}_{0}.$\textit{The
following conditions are equivalent for a proper cut} $\mathrm{I}$ \textit{of%
} $\mathcal{M}$:\medskip

\noindent \textbf{(1)} \textit{There is a self-embedding} $j$ \textit{of} $(%
\mathcal{M},\mathcal{A})$ \textit{such that }$\mathrm{Fix}(j)=\mathrm{I}.$%
\medskip

\noindent \textbf{(2) }$\mathrm{I}$\ \textit{is a strong cut of }$\mathcal{M}
$\textit{\ and }$\mathrm{I}\prec _{\Sigma _{1}}\mathcal{M}.$\medskip

\noindent \textbf{(3) }\textit{There is a proper initial self-embedding} $j$
\textit{of} $(\mathcal{M},\mathcal{A})$ \textit{such that }${\mathrm{Fix}(j)=%
\mathrm{I}.}$\medskip

Theorem A is established in Section 3, and Section 4 is devoted to proving
Theorem B.\medskip

\noindent \textbf{Acknowledgment.}~I am indebted to my PhD supervisor, Ali
Enayat, for his encouragements, precious comments and feedback in the
formation of this paper.

\section{Preliminaries}

In this section we review some definitions and basic results that are
relevant to the statements and proofs of our main results.

\begin{itemize}
\item $\mathrm{WKL}_{0}$ is the second order theory whose models are of the
form of $(\mathcal{M},\mathcal{A}),$ where $(\mathcal{M},\mathcal{A})$
satisfies (1) Induction for $\Sigma _{1}^{0}$ formulas; (2) Comprehension
for $\Delta _{1}^{0}$-formulas; and (3) Weak K\"{o}nig's Lemma (which
asserts that every infinite subtree of the full binary tree has an infinite
branch).

\item $\mathrm{Exp}:=\forall x\exists y\ \mathrm{Exp}(x,y)$, where $\mathrm{%
Exp}(x,y)$ is a $\Delta _{0}$-formula that expresses $2^{x}=y$ within
\textrm{I}$\Delta _{0}.$

\item The binary $\Delta _{0}$-formula $xEy$, known as \textit{Ackermann's
membership relation}, expresses \textquotedblleft the $x$-th bit of the
binary expansion of $y$ is 1\textquotedblright\ within \textrm{I}$\Delta
_{0}.$

\item A subset $\mathrm{X}$ of $\mathrm{M}$ is \textit{coded} in $\mathcal{M}
$ iff there is some $a\in \mathrm{M}$ such that $\mathrm{X}={\left(
a_{E}\right) ^{\mathcal{M}}:=\{x\in \mathrm{M}:\mathcal{M}\models xEa\}}.$
\newline
Given $a\in \mathrm{M},$ by $\underline{a}$ we mean the set $\{x\in \mathrm{M%
}:x<a\}$. Note that $\underline{a}$ is coded in $ \mathcal{M} $, where  $\mathcal{M}$ is a model of $ \mathrm{I%
}\Delta _{0}$, provided $2^{a}$ exists in $\mathcal{M}$.\newline
It is well-known that for any $n>0$ and $\mathcal{M}\models \mathrm{I}\Sigma
_{n}$, if $\varphi (x,a)$ is an unary $\Sigma _{n}$-formula, where $a$ is a
parameter from $\mathcal{M}$, then the set ${\varphi ^{\mathcal{M}}(a):={\{}}%
b\in \mathrm{M}{{:\mathcal{M}\models \varphi (b,a)\}}}$ is \textit{%
piece-wise coded} in $\mathcal{M}$; i.e. for every $c\in \mathrm{M}$, $%
\varphi ^{\mathcal{M}}(a)\cap \underline{c}$ is coded in $\mathcal{M}$. More
specifically, there is some element less than $2^{c}$ which codes $\varphi ^{%
\mathcal{M}}(a)\cap \underline{c}$. \newline
Moreover, the above statement holds for $n=0$ if $\mathcal{M}\models \mathrm{%
I}\Delta _{0}+\mathrm{Exp}$.

\item For every cut $\mathrm{I}$ of $\mathcal{M}$, the $\mathrm{I}$\textit{%
-standard system} of $\mathcal{M}$, presented by $\mathrm{SSy}_{I}(\mathcal{M%
})$, is the family consisting of sets of the form $a_{E}\cap \mathrm{I}$,
where $a\in \mathrm{M}$; in other words:
\begin{equation*}
\mathrm{SSy}_{\mathrm{I}}(\mathcal{M})=\{\left( a_{E}\right) ^{\mathcal{M}%
}\cap \mathrm{I}:\ a\in \mathrm{M}\}.
\end{equation*}%
When $\mathrm{I}$ is the standard cut, i.e. $\mathrm{I}=\mathbb{N}$, we
simply write $\mathrm{SSy}(\mathcal{M})$ instead of $\mathrm{SSy}_{\mathbb{N}%
}(\mathcal{M}).$

\item Given a proper cut $ \mathrm{I} $ of $ \mathcal{M} $, $ \mathrm{I} $ is called a \textit{strong cut}, if for every coded function $ f $ in $ \mathcal{M} $ whose domain contains $ \mathrm{I} $, there exists some $ s\in \mathrm{M} $ such that for every $ i\in \mathrm{I} $ it holds that $ f(i)\notin \mathrm{I} $ iff $ s<f(i) $.

\item $\mathrm{Sat}_{_{\Sigma _{n}}}$ is the arithmetical formula defining
the satisfaction predicate for $\Sigma _{n}$-formulas within  $\mathrm{I}%
\Delta _{0}+\mathrm{Exp}$. It is well-known that for each positive $n\in
\omega $, $\mathrm{Sat}_{_{\Sigma _{n}}}$ can be expressed by a $\Sigma _{n}$%
-formula in $\mathrm{I}\Sigma _{n}$; furthermore, within a model of  ${\mathrm{I}%
\Delta _{0}+\mathrm{Exp}}$ (with the help of a nonstandard parameter if the model is nonstandard), $\mathrm{Sat}_{_{\Sigma _{0}}}$ (which is also written as $\mathrm{%
Sat}_{_{\Delta _{0}}}$) is expressible both as a $\Sigma _{1}$ and $\Pi _{1}$%
-formula.

\item The strong $\Sigma _{n}$-\textit{Collection scheme }consists of
formulas of the following form where $\varphi $ is a $\Sigma _{n}$-formula:

$\forall w\ \forall v\ \exists z\forall x<v(\exists y\ \varphi
(x,y,w)\rightarrow \exists y<z\ \varphi (x,y,w)).$

It is well-known that the strong $\Sigma _{n}$-Collection scheme is provable
in $\mathrm{I}\Sigma _{n}$ for every $n>0.$

\item Every model $\mathcal{M}$ of $\mathrm{I}\Delta _{0}+\mathrm{Exp}$
satisfies the \textit{Coded Pigeonhole Principle}, i.e. if $b\in \mathrm{M}$%
, and $f:\underline{b+1}\rightarrow \underline{b}$ is a coded function in $%
\mathcal{M}$, then $f$ is not injective.

\item By an embedding $j$ from second order model $(\mathcal{M},\mathcal{A})$
into $(\mathcal{N},\mathcal{B})$, we mean that $j$ is an embedding from $%
\mathcal{M}$ into $\mathcal{N}$ such that for every $X\subseteq M$, \ $X\in
\mathcal{A}$ iff $j(X)=Y\cap j(M)$ for some $Y\in \mathcal{B}$.

\item When $\mathcal{M}$ and $\mathcal{N}$ are models of arithmetic and $b$
is in $\mathrm{M}$, we write $\mathcal{M}\subseteq _{\mathrm{end,\Pi }%
_{1,\leq b}}\mathcal{N}$, when $\mathcal{N}$ is an end extension of $%
\mathcal{M}$ and all $\Pi _{1}$-formulas whose parameters are in $\underline{%
b+1}$ are absolute in the passage between $\mathcal{M}$ and $\mathcal{N}$,
i.e., $\mathrm{Th}_{\Pi _{1}}(\mathcal{M},m)_{m\leq b}=\mathrm{Th}_{\Pi
_{1}}(\mathcal{N},m)_{m\leq b}.$
\end{itemize}

\noindent \textbf{Theorem 2.1.}~\textbf{(\cite[Theorem 3.2]{Me Tanaka})} \textit{Let }%
$(\mathcal{M},\mathcal{A})$\textit{\ be a countable model of }$\mathrm{WKL}%
_{0}$\textit{\ and let }$b\in \mathrm{M}$\textit{. Then }$\mathcal{M}$%
\textit{\ has a countable recursively saturated proper end extension }$\mathcal{N}$%
\textit{\ satisfying }$\mathrm{I}\Delta _{0}+\mathrm{Exp}+\mathrm{B}\Sigma
_{1}$ \textit{such that }$\mathrm{SSy}_{\mathrm{M}}(\mathcal{N})=\mathcal{A}$%
\textit{, and }$\mathcal{M}\subseteq _{\mathrm{end,\Pi_{1},\leq b}}\mathcal{%
N}$.\textit{\medskip }

\noindent \textbf{Remark 2.2.}~The proofs of our main results take advantage
of the following additional features of the model $\mathcal{N}$ constructed
in Enayat's \textit{proof} of Theorem 2.1, namely: given $b$ in $\mathcal{M}$%
, there is an elementary chain of models ${(\mathcal{N}_{n}:n\in \omega )}$
satisfying the following three properties:

\begin{enumerate}
\item[$(i)$] $\mathcal{N}=\cup _{n\in \omega }\mathcal{N}_{n}$;

\item[$(ii)$] For every $n\in \omega $,  $\mathcal{M}\subseteq _{\mathrm{end,\Pi_{1},\leq b}}\mathcal{N%
}_{n}\prec \mathcal{N}$;

\item[$(iii)$] For every $n\in \omega $ the elementary diagram of $(\mathcal{%
N}_{n},a)_{a\in \mathrm{N}_{n}}$ is available in $(\mathcal{M},\mathcal{A})$
via some $\mathrm{ED}_{n}\in \mathcal{A}$. Note that $\mathrm{Th}\left( (%
\mathcal{N}_{n},a)_{a\in \mathrm{N}_{n}}\right) $ is a proper subset of $%
\mathrm{ED}_{n}$ since $\mathrm{ED}_{n}$ includes sentences of nonstandard
length.$\medskip $
\end{enumerate}

\noindent \textbf{Remark 2.3. }Enayat \cite{Me Tanaka} noted that if $(%
\mathcal{M},\mathcal{A})$\ is a model of $\mathrm{WKL}_{0}$ and $b$ is in $%
\mathrm{M}$, and there is some end extension $\mathcal{N}$ of $\mathcal{M}$
such that (1) $\mathcal{N}\models I\Delta _{0}+\mathrm{Exp}$, (2) $\mathrm{%
SSy}_{\mathrm{M}}(\mathcal{N})=\mathcal{A}$, and (3) there is an initial
self-embedding $j_{1}\ $of $\mathcal{N}$ onto an initial segment that is
bounded above by $b$, then the restriction $j$ of $j_{1}$ to $\mathcal{M}$
is an embedding of $\mathcal{M}$ onto an an initial segment $\mathrm{J}$ of $%
\mathcal{M}$ that is below $b$ which has the important feature that $j$ is
an isomorphism between $(\mathcal{M},\mathcal{A})$ and $(\mathrm{J},\mathcal{%
A}_{\mathrm{J}})$. Note that if $\mathrm{I}$ is a proper cut of $\mathcal{M}%
, $ then $\mathrm{I}_{\mathrm{fix}}(j_{1})=\mathrm{I}$ implies that $\mathrm{%
I}_{\mathrm{fix}}(j)=\mathrm{I}$; and $\mathrm{Fix}(j_{1})=\mathrm{I}$\
implies that $\mathrm{Fix}(j)=\mathrm{I}.$

\begin{itemize}
\item The following theorem summarizes some of the results about $\mathrm{I}%
_{\mathrm{fix}}(j)$ and $\mathrm{Fix}(j)$ from \cite{our} which will be
employed in this paper:
\end{itemize}

\noindent \textbf{Theorem 2.4.}~Suppose $\mathcal{M}\models \mathrm{I}\Delta
_{0}+\mathrm{Exp}$ \textit{and} $j$ \textit{is a nontrivial self-embedding of%
} $\mathcal{M}$. \textit{Then:}

\noindent \textbf{(a)} $\mathrm{I}_{\mathrm{fix}}(j)\models \mathrm{I}\Delta
_{0}+\mathrm{B}\Sigma _{1}+\mathrm{Exp}$.

\noindent \textbf{(b)} \textit{If} $\mathcal{M}\models \mathrm{I}\Sigma _{1}$
\textit{then} $\mathrm{Fix}(j)$ \textit{is a} $\Sigma _{1}$-\textit{%
elementary submodel of} $\mathcal{M}$. \textit{Moreover, if} $\mathrm{Fix}%
(j) $ \textit{is a proper initial segment of} $\mathcal{M}$, \textit{then it
is a strong cut} of $\mathcal{M}$.

\begin{itemize}
\item Given two countable nonstandard model $\mathcal{M}$ and $\mathcal{N}$
of $\mathrm{I}\Sigma _{1}$ which share a common proper cut $\mathrm{I}$, the
following theorem from \cite[Cor. 3.3.1]{our} provides a useful sufficient
condition for existence of a proper initial embedding between $\mathcal{M}$
and $\mathcal{N}$ which fixes each element of $\mathrm{I}$:\medskip
\end{itemize}

\noindent \textbf{Theorem 2.5.}~\textit{Let} $\mathcal{M}$, $\mathcal{N}$
\textit{and} $I$ \textit{be as above such that }$I$\textit{\ is closed under
exponentiation. The following are equivalent}:\textit{\smallskip }

\noindent \textbf{(1)}\textit{\ There is a proper initial embedding }$f$%
\textit{\ of }$\mathcal{M}$ \textit{into} $\mathcal{N}$ \textit{such that} $%
f(i)=i$ \emph{for all }$i\in I.$\textit{\smallskip }

\noindent \textbf{(2) }$\mathrm{Th}_{\Sigma _{1}}(\mathcal{M},i)_{i\in
I}\subseteq \mathrm{Th}_{\Sigma _{1}}(\mathcal{N},i)_{i\in I}$ \textit{and} $%
\mathrm{SSy}_{I}(\mathcal{M})=\mathrm{SSy}_{I}(\mathcal{N}).$\textit{%
\medskip }

\begin{itemize}
\item Another prominent subsystem of second order arithmetic is $\mathrm{%
ACA_{0}}$, in which the comprehension scheme is restricted to formulas with
no second order quantifier. The following results of Paris and Kirby \cite%
{Jeff and Laurie} and Gaifman \cite[Thm. 4.9-4.11]{Gaifman} concerning $%
\mathrm{ACA_{0}}$ are employed in the proof of Theorem B.\footnote{%
Gaifman couched his results in terms of arbitrary models of $\mathrm{PA}(%
\mathcal{L})$\ for countable $\mathcal{L}$. Note that if $(\mathcal{M},%
\mathcal{A})\models \mathrm{ACA}_{0}$, then the expansion $(\mathcal{M}%
,A)_{A\in \mathcal{A}}$ of $\mathcal{M}$ is a model of $\mathrm{PA}(\mathcal{%
L})$, where $\mathcal{L}$ is the extension of $\mathcal{L}_{A}$ by predicate
symbols for each $A\in \mathcal{A}.$ Moreover, the collection of subsets of $%
\mathrm{M}$ that are parametrically definable in $(\mathcal{M},A)_{A\in
\mathcal{A}}$ coincides with $\mathcal{A}$.}
\end{itemize}

\noindent \textbf{Theorem 2.6.}~\textbf{(Paris and Kirby)} Suppose $\mathcal{M}%
\models \mathrm{I}\Delta _{0}$. \textit{The following are equivalent for a
proper cut} $\mathrm{I}$ \textit{of }$\mathcal{M}$:\textit{\smallskip }

\noindent \textbf{(a)} $\mathrm{I}$ \textit{is a strong cut of} $\mathcal{M}$%
.\textit{\smallskip }

\noindent \textbf{(b)} $(\mathrm{I},\mathrm{SSy}_{\mathrm{I}}(\mathcal{M}%
))\models \mathrm{ACA}_{0}.$\medskip

\noindent \textbf{Theorem 2.7.}~\textbf{(Gaifman)} \textit{Given a
countable model }$(\mathcal{M},\mathcal{A})$\textit{\ of }$\mathrm{ACA}_{0}$
\textit{and a linear order }$\mathbb{L}$, \textit{there exists an end
extension} $\mathcal{M}_{\mathbb{L}}$ \textit{of} $\mathcal{M}$ \textit{such
that there is an isomorphic copy} $\mathbb{L}^{\prime }=\{c_{l}:l\in \mathbb{%
L\}}$ \textit{of} $\mathbb{L}$ \textit{in} $\mathrm{M}_{\mathbb{L}%
}\backslash \mathrm{M}$, \textit{and there is a composition preserving
embedding }$j\mapsto \widehat{j}$ \textit{from the semi-group of initial
self-embeddings of} $\mathbb{L}$ \textit{into the semi-group of initial
self-embeddings of }$\mathcal{M}_{\mathbb{L}}$ \textit{that satisfy the
following properties}:\textit{\smallskip }

\noindent \textbf{(a) }$\mathrm{SSy}_{\mathrm{M}}(\mathcal{M}_{\mathbb{L}})=%
\mathcal{A}$ \textit{and }$\mathrm{Fix}(\widehat{j})=M$ \textit{for each
initial self-embedding} $j$ \textit{of} ${\mathbb{L}}$ \textit{that is fixed
point free}.\textit{\smallskip }

\noindent \textbf{(b) }\textit{For each initial self-embedding} $j$ \textit{%
of} ${\mathbb{L}}$, $\widehat{j}$\textit{\ is an \textbf{elementary} initial
self-embedding} \textit{of} $\mathcal{N}_{\mathbb{L}}$, \textit{i.e.} $%
\widehat{j}(\mathcal{M}_{\mathbb{L}})\preceq _{\mathrm{end}}\mathcal{N}_{%
\mathbb{L}}$.\textit{\smallskip }

\noindent \textbf{(c) }$\mathbb{L}^{\prime }$ \textit{is downward cofinal in
}$M_{\mathbb{L}}\backslash M$ \textit{if} $\mathbb{L}$ \textit{has no first
element. \smallskip }

\noindent \textbf{(d)} \textit{For any} $l_{0}\in \mathbb{L}$, $l_{0}$
\textit{is a strict upper bound for} $j(\mathbb{L})$ \textit{iff} $c_{l_{0}}$
\textit{is a strict upper bound for} $\widehat{j}(M_{\mathbb{L}})$.\smallskip

\section{The longest cut fixed by self-embeddings}

This section is devoted to the proof of Theorem A. Recall that if $(\mathcal{%
M},\mathcal{A})$ is a model of $\mathrm{WKL}_{0}$ there are arbitrary large
as well as arbitrary small nonstandard cuts in $\mathcal{M}$ that are closed
under exponentiation. More specifically, for every nonstandard $a$ in $%
\mathrm{M}$, there are nonstandard cuts $\mathrm{I}_{1}$ and $\mathrm{I}_{2}$ (as defined
below) such that $\mathrm{I}_{1}<a\in \mathrm{I}_{2}$, and both are closed
under exponentiation:

\begin{center}
$\mathrm{I}_{1}:=\{x\in M:2_{n}^{x}<a\ \text{for all}\ n\in \omega \}$,
where $2_{0}^{x}:=x$, and for every $n\in \omega $, $%
2_{n+1}^{x}:=2^{2_{n}^{x}}$;\medskip \\[0pt]
$\mathrm{I}_{2}:=\{x\in M:x<2_{n}^{a}\ \text{for some}\ n\in \omega \}$.
\end{center}

\noindent So Theorem A implies that $\mathrm{I_{fix}}(j)$ can be arranged to
be as high or as low in the nonstandard part of $\mathcal{M}$ as desired. In
particular, Theorem A is a strengthening of Tanaka's Theorem.\newline
\medskip

\noindent \textbf{Proof of Theorem A.}~$(1)\Rightarrow (2)$ is an immediate
consequence of Theorem 2.4.(a), and $(3)\Rightarrow (1)$ is trivial so we
concentrate on establishing $(2)\Rightarrow (3)$. \medskip

Assume that $\mathrm{I}$ is closed under exponentiation and fix some $a\in
\mathrm{M}\setminus\mathrm{I}$. We leave it as an exercise for the reader to
use strong $\Sigma _{1}$-Collection along with the fact that \textrm{Sat}$%
_{\Delta _{0}}$ has a $\Sigma _{1}$-description in $\mathcal{M}$ to show
that there is some $b\in \mathrm{M}$ such that:\medskip

\noindent $(\sharp )\ \ \ $ $\mathcal{M}\models \forall w<a\left( \exists z\
\delta (z,w)\rightarrow \exists z<b\ \delta (z,w\right)) $, for all $\Delta
_{0}$-formulas $\delta $.\medskip

\noindent Next we invoke Theorem 2.1 to get hold of a countable recursively
saturated proper end extension $\mathcal{N}$ of $\mathcal{M}$ such that $%
\mathcal{N}\models \mathrm{I}\Delta _{0}+\mathrm{B}\Sigma _{1}+\mathrm{Exp}$%
, ${\mathrm{SSy}_{\mathrm{M}}(\mathcal{N})=\mathcal{A}}$, and $\mathcal{M}%
\subseteq _{\mathrm{end,\Pi_{1},\leq b}}\mathcal{N}$. Moreover, we will
safely assume that the model $\mathcal{N}$ additionally satisfies the three
properties listed in Remark 2.2. In light of Remark 2.3, in order to
establish (3) it suffices to construct a proper initial self-embedding $j$
of $\mathcal{N}$ such that $j(\mathrm{N})<b$ and $\mathrm{I}_{\mathrm{fix}}(j)=%
\mathrm{I}$. The construction of the desired $j$ is the novel element of the
proof of Theorem A, which we now turn to.\medskip

To construct $j$ we will employ a modification of the strategy employed in
the proof of $(2)\Rightarrow (3)$ of \cite[Theorem 4.1]{our}, using a
3-level back-and-forth method. A modification is needed since we need to
overcome the fact that I$\Sigma _{1}$ need not hold in $\mathcal{N}$;
instead we will rely on recursive saturation of $\mathcal{N}$ and the
properties of $\mathcal{N}$ listed in Remark 2.2. First, note that $(\sharp
) $ together with the fact that $\mathrm{Th}_{\Pi _{1}}(\mathcal{M}%
,x)_{x\leq b}=\mathrm{Th}_{\Pi _{1}}(\mathcal{N},x)_{x\leq b}$,
implies:\medskip

\noindent $(\ast )\ \ \ $ $\mathcal{N}\models \forall w<a\left( \exists z\
\delta (z,w)\rightarrow \exists z<b\ \delta (z,w\right)) $, for all $\Delta
_{0}$-formula $\delta $.\medskip

\noindent Since $\mathrm{I}$\ is closed under exponentiation, we can choose $%
\{c_{n}:n\in \omega \}$ that is downward cofinal in $\mathrm{M}\setminus
\mathrm{I}$ such that $c_{0}=a$ and $2^{c_{n+1}}<c_{n}$ for all $n\in \omega
$.\medskip

The proof will be complete once we recursively construct finite sequences $%
\bar{u}:=\left( u_{0},...,u_{m-1}\right) $ and $\bar{v}:=\left(
v_{0},...,v_{m-1}\right) $ of elements of $\mathcal{N}$ for all $n\in \omega
$ such that:\medskip

\begin{itemize}
\item[$(i)$] $u_{0}=0=v_{0}$.

\item[$(ii)$] For every $c$ in $\mathrm{N}$ there is some $n\in \omega $
such that $c=u_{n}$.

\item[$(iii)$] For every $n\in \omega $, $v_{n}<b$, and if for some $c$ in $%
\mathrm{N}$ it holds that $c<v_{n}$, then there is some $m\in \omega $ such that $c=v_{m}
$.

\item[$(iv)$] For every $m\in \omega $ the following condition holds:

$(\ast _{m}):\ \ \ $ $\mathcal{N}\models \forall w<c_{m}\ \left( \exists z\
\delta (z,w,\bar{u})\rightarrow \exists z<b\ \delta (z,w,\bar{v})\right) $,%
\newline
for every $\Delta _{0}$-formula $\delta $.

\item[$(v)$] For every $m\in \omega $, there is some $n\in \omega $ such
that $u_{n}<c_{m}$ and $u_{n}\neq v_{n}$.
\end{itemize}

Note that $(\ast _{0})$ holds thanks to $(\ast )$ since $c_{0}=a.$ Let $%
\{a_{n}:n\in \omega \}$ and $\{b_{n}:n\in \omega \}$ respectively be
enumerations of element of $\mathrm{N}$ and $\underline{b}$. By statement $(i)$ and $%
(\ast _{0})$ the first step of induction holds. Suppose for $m\in \omega $, $%
\bar{u}$ and $\bar{v}$ are constructed such that $(\ast _{m})$ holds. In
order to find suitable $u_{m+1}$ and $v_{m+1}$, by considering congruence
modulo 3 we have three cases for $m+1$: Case 0 takes care of $(ii)$ and $%
(iv) $, Case 1 takes care of $(iii)$ and $(iv)$, and Case 2 takes care of $%
(v)$ and $ (iv) $.\medskip

\textbf{CASE 0 ($\mathbf{m+1=3k}$, for some $\mathbf{k}\in \omega $):} In
this case if $a_{k}$ is one of the elements of $\bar{u}$, put $u_{m+1}=u_{m}$
and $v_{m+1}=v_{m}$. Otherwise, put $u_{m+1}=a_{k}$ and define:

\begin{center}
$p(y):=\{y<b\}\cup \left\{ \forall w<c_{m+1}(\exists z\ \delta (z,w,\bar{u}%
,a_{k})\rightarrow \exists z<b\ \delta (z,w,\bar{v},y)):\delta \ \mathrm{is\
a}\ \Delta _{0}\mathrm{-formula}\right\} .$
\end{center}

\noindent Note that $p(y)$ is a recursive type. Since $\mathcal{N}$ is
recursively saturated, it suffices to prove that $p(y)$ is finitely
satisfiable and let $v_{m+1}$ be one of the realizations of $p(y)$ in $%
\mathcal{N}$. Since $p(y)$ is closed under conjunctions we only need to show
that each formula in $p(y)$ is satisfiable. For this purpose, suppose $%
\delta $ is a $\Delta _{0}$-formula, and let

\begin{center}
$\mathrm{D}:=\left\{ w\in \underline{c_{m+1}}:\mathcal{N}\models \exists z\
\delta (z,w,\bar{u},a_{k})\right\} .$
\end{center}

\noindent We claim that there is some $d<2^{c_{m+1}}$ which codes $\mathrm{D}
$ in $\mathcal{N}$. To see this, we note that in the above definition, $%
\mathcal{N} $ can be safely replaced by some $\mathcal{N}_{n}$, where $n$ is
large enough to contain the parameters $\bar{u}$ and $a_{k}$ (thanks to
properties $(i)$ and $(ii)$ in Remark 2.2). On the other hand, by property $%
(iii)$ in Remark 2.2, there is some $\mathrm{ED}_{n}\in \mathcal{A}$ such
that:

\begin{center}
$\mathrm{D}=\left\{ w\in \underline{c_{m+1}}:\ulcorner \exists z\ \delta
(z,w,\bar{u},a_{k})\urcorner \in \mathrm{ED}_{n}\right\} .$
\end{center}

\noindent Since $(\mathcal{M},\mathcal{A})$ satisfies $\mathrm{I}\Sigma
_{1}^{0}$, the above characterization of $\mathrm{D}$ shows that $\mathrm{D}$
is coded in $\mathcal{M}$ (and therefore in $ \mathcal{N} $) by some $d<2^{c_{m+1}}$ (recall that the code of
each subset of $\underline{m}$ is below $2^{m}).$ Therefore we have:\medskip

\noindent (1) $\mathcal{N}\models \forall w<c_{m+1}\left( wEd\rightarrow
\exists z\ \delta (z,w,\bar{u},a_{k})\right) ;$\medskip

\noindent By putting (1) together with $\mathrm{B}\Sigma _{1}$ in $\mathcal{N%
}$, and existentially quantifying $a_{k}$ we obtain:\medskip

\noindent $(2)$ $\mathcal{N}\models \exists t,x\ \forall w<c_{m+1}\left(
wEd\rightarrow \exists z<t\ \delta (z,w,\bar{u},x)\right) .$\medskip

\noindent On the other hand, coupling (2) with $(\ast _{m})$ yields:\medskip

\noindent $(3)$ $\mathcal{N}\models \exists t,x<b \ \forall
w<c_{m+1}((wEd\rightarrow \exists z<t\ \delta (z,w,\bar{v},x))),$\medskip

\noindent which makes it clear that each formula in $p(y)$ is satisfiable in
$\mathcal{N}$.\medskip

\textbf{CASE 1 ($\mathbf{m+1=3k+1}$, for some $\mathbf{k}\in \omega $):} In
this case if ${b_{k}\geq \mathrm{Max}\{\bar{v}\}}$ or if it is one of the
elements of $\bar{v}$, put $u_{m+1}=u_{m}$ and $v_{m+1}=v_{m}$. Otherwise,
put $v_{m+1}=b_{k}$ and define:

\begin{center}
$q(x):=\left\{ \forall w<c_{m+1}(\forall z<b\ \lnot \delta (z,w,\bar{v}%
,b_{k})\rightarrow \forall z\ \lnot \delta (z,w,\bar{u},x)):\delta \ \mathrm{%
is\ a}\ \Delta _{0}\mathrm{-formula}\right\} .$
\end{center}

\noindent $q(x)$ is clearly a recursive type and closed under conjunctions,
so by recursive saturation of $\mathcal{N}$ it suffices to verify that each
formula in $q(x)$ is satisfiable in $\mathcal{N}$, and let $u_{m+1}$ be one
of the realizations of $q(x)$ in $\mathcal{N}$. Suppose some formula in $q$
is not realizable in $\mathcal{N}$, then for some $\Delta _{0}$-formula $%
\delta $ we have:

\begin{center}
$\mathcal{N}\models \forall x\left( \exists w<c_{m+1}\left( \forall z<b\
\lnot \delta (z,w,\bar{v},b_{k})\wedge \exists z\ \delta (z,w,\bar{u}%
,x)\right) \right) .$
\end{center}

\noindent Let:

\begin{center}
$\mathrm{R}:=\left\{ w\in \underline{c_{m+1}}:\mathcal{N}\models \forall
z<b\ \lnot \delta (z,w,\bar{v},b_{k})\right\} .$
\end{center}

\noindent Since $\mathrm{R}$ is $\Delta _{0}$-definable in $\mathcal{N}$
there exists some $r<2^{c_{m+1}}$ which codes $\mathrm{R}$ in $\mathcal{N%
}$. Therefore, \medskip

\noindent (4) $\mathcal{N}\models \forall x<\mathrm{\mathrm{Max}}\{\bar{u}%
\}\left( \exists w<c_{m+1}\left( wEr\wedge \exists z\ \delta (z,w,\bar{u}%
,x)\right) \right) ,$\medskip

\noindent which by $\Sigma _{1}$-Collection in $\mathcal{N}$ implies:\medskip

\noindent (5) $\mathcal{N}\models \exists t\ \forall x<\mathrm{\mathrm{Max}}%
\{\bar{u}\}\left( \exists w<c_{m+1}\left( wEr\wedge \exists z<t\ \delta (z,w,%
\bar{u},x)\right) \right) .$\medskip

\noindent Putting (5) together with $(\ast_{m})$ yields:\medskip

\noindent $(6)$ $\mathcal{N}\models \exists t<b\ \forall x<\mathrm{\mathrm{%
Max}}\{\bar{v}\}\left( \exists w<c_{m+1}\left( wEr\wedge \exists z<t\ \delta
(z,w,\bar{v},x)\right) \right) .$\medskip

\noindent By substituting $b_{k}$ for $x$ in (6) we obtain:\medskip

\noindent (7) $\mathcal{N}\models \exists t<b\ \forall x<\mathrm{\mathrm{Max}%
}\{\bar{v}\}\left( \exists w<c_{m+1}\left( wEr\wedge \exists z<t\ \delta
(z,w,\bar{v},b_{k})\right) \right) .$\medskip

\noindent But (7) contradicts the assumption that $r$ codes $\mathrm{R}$. So
$q(x)$ is finitely satisfiable.\medskip

\textbf{CASE 2 ($\mathbf{m+1=3k+2}$, for some $\mathbf{k}\in \omega $):}
Consider the type $l(x,y):=\left\{ {x\neq y,x\leq c_{k}}\right\}\cup
l_{0}(x,y)$, where:

\begin{center}
$l_{0}(x,y):=\left\{ {\forall w<c_{m+1}}\left( {\exists z\ \delta (z,w,\bar{u%
},x)\rightarrow \exists z<b\ \delta (z,w,\bar{v},y)}\right) :\delta \
\mathrm{is\ a}\ \Delta _{0}\mathrm{-formula}\right\} .$
\end{center}

\noindent Once we demonstrate that $l(x,y)$ is realized in $\mathcal{N}$ we
can define $(u_{m+1},v_{m+1})$ as any realization in $\mathcal{N}$ of ${%
l(x,y)}$. Since $l_{0}(x,y)$ is closed under conjunctions and $\mathcal{N}$
is recursively saturated, to show that $l(x,y)$ is realized in $\mathcal{N}$
it suffices to demonstrate that the conjunction of ${x\neq y}$ and ${x\leq
c_{k}}$, and each formula in $l_{0}(x,y)$ is satisfiable in $\mathcal{N}.$
So suppose $\delta $ is a $\Delta _{0}$-formula and for each $s<c_{k}$
consider the map $F$ from $\underline{c_{k}}$ to the power set of $%
\underline{c_{m+1}}$ by:

\begin{center}
$F(s):={\{w\in \underline{c_{m+1}}:\ \mathcal{N}\models \exists z\ }${$\delta $}${(z,w,\bar{u},s)\}}.$\\[0pt]
\end{center}

\noindent Thanks to properties $(i)$ through $(iii)$\ of $\mathcal{N}$
listed in Remark 2.2, there is some $\mathrm{ED}_{n}\in \mathcal{A}$ such
that:

\begin{center}
$F(s)={\{w\in \underline{c_{m+1}}:\ulcorner  \exists z\ }${%
$\delta $}${(z,w,\bar{u},s)\urcorner }\in \mathrm{ED}_{n}{\}}.$\\[0pt]
\end{center}

\noindent Since $(\mathcal{M},\mathcal{A})$ satisfies $\mathrm{I}\Sigma
_{1}^{0}$, the above characterization of $F(s)$ together with the
veracity of $\mathrm{I}\Sigma _{1}^{0}$ in $(\mathcal{M},\mathcal{A})$ makes
it clear that $F$ is coded in $\mathcal{M}$ by some $f$ (and therefore in $%
\mathcal{N}$) that codes a function from ${\underline{{c_{k}}}}$ to $%
\underline{2^{c_{m+1}}}$ with ${f(s):=\sum_{l\in F(s)}2^{l}}.$ On the other
hand the definition of $f(s)$ and the assumption that $2^{c_{n+1}}<c_{n}$
for all $n\in \omega $ makes it clear that:

\begin{center}
$f(s)\leq \sum_{l<c_{m+1}}2^{l}=2^{c_{m+1}}-1<2^{c_{m+1}}<c_{m}<c_{k}$.
\end{center}

\noindent So by the coded pigeonhole principle there are \textbf{distinct} ${%
s,s^{\prime }<c_{k}}$ such that $f(s)=f(s^{\prime })$, in other words:

\begin{center}
$\mathcal{N}\models \forall w<c_{m+1}(\exists z\ \delta (z,w,\bar{u}%
,s)\leftrightarrow \exists z\ \delta (z,w,\bar{u},s^{\prime })).$
\end{center}

\noindent Now by repeating the argument used in Case 0 for $\left( \bar{u}%
,s,s^{\prime }\right) $ we can find some $t,t^{\prime }<b$ such that:

\begin{center}
$\mathcal{N}\models \forall w<c_{m+1}(\exists z\ \delta (z,w,\bar{u}%
,s,s^{\prime })\rightarrow \exists z<b\ \delta (z,w,\bar{v},t,t^{\prime })).$
\end{center}

\noindent Since $s\neq s^{\prime }$, either $s\neq t$ or $s\neq t^{\prime }$%
. So the conjunction of ${x\neq y,}$ ${x\leq c_{k}}$, and each formula in $%
l_{0}(x,y)$ is satisfiable in $\mathcal{N}$, and the proof is
complete.\hfill $\square $\medskip

\section{Cuts which are fixed-point sets of self-embeddings}

In this section we present the proof of Theorem B. But before going through
the proof, let us point out that a model of $\mathrm{WKL}_{0}$ does not
necessarily carry a cut satisfying statement (2) of Theorem B (see \cite[%
Remark 5.1.1]{our} for an explanation). However, if $\mathcal{M}\models
\mathrm{PA}$, there are arbitrarily high strong cuts $\mathrm{I}$ in $%
\mathcal{M}$ such that $\mathrm{I}\prec _{\Sigma _{1}}\mathcal{M}$. To see
this when $\mathcal{M}$ is a countable model of $\mathrm{PA}$, let $\mathcal{%
A}$ be the family of definable subsets of $\mathcal{M}$. Since $(\mathcal{M},%
\mathcal{A})\models \mathrm{ACA}_{0}$ and $\mathrm{WKL}_{0}$ is a subsystem
of $\mathrm{ACA}_{0}$, by Theorem 1.1 (Tanaka's theorem), for every $a\in
\mathrm{M}$ there is a cut $\mathrm{I}$ containing $a$ such that $(\mathcal{M%
},\mathcal{A})\cong (\mathrm{I},\mathcal{A}_{\mathrm{I}})$ and $\mathrm{%
I\prec _{\Sigma _{1}}}\mathcal{M}$. Furthermore, $\mathrm{I}$ is strong cut
of $\mathcal{M}$ by Theorem 2.6 since $\mathcal{A}_{\mathrm{I}}=\mathrm{SSy}%
_{\mathrm{I}}(\mathcal{M})$.\medskip

\noindent \textbf{Proof of Theorem B.}~$(1)\Rightarrow (2)$ is an immediate
consequence of Theorem 2.4.(b), and $(3)\Rightarrow (1)$ is trivial; so we
concentrate on the proof of $(2)\Rightarrow (3)$. \medskip

\noindent Suppose $\mathrm{I}$ is a strong cut in $\mathcal{M}$ and $\mathrm{%
I}\prec _{\mathrm{end},\ \Sigma _{1}}\mathcal{M}$. The proof of (3) is
inspired by the proof of \cite[Theorem 5.1]{our} and consists of the
following four stages:\medskip

\textbf{Stage 1:} Fix some $b_{0}\in \mathrm{M}\setminus \mathrm{I}$. Using
Theorem 2.1, let $\mathcal{N}$ be a model of $\mathrm{I}\Delta _{0}+\mathrm{B%
}\Sigma _{1}+\mathrm{Exp}$ such that $\mathrm{SSy}_{\mathrm{M}}(\mathcal{N})=%
\mathcal{A}$, $\mathcal{M}\subseteq _{\mathrm{end,\Pi_{1},\leq b_{0}}}%
\mathcal{N}$, and the three conditions specified in Remark 2.2 hold for $%
\mathcal{N}$. \medskip

\textbf{Stage 2:} Let $\mathbb{Q}$ be the set of rational numbers with its
natural ordering. Since $\mathrm{I}$ is a strong cut in $\mathcal{M}$, by
Theorem 2.6, and the case $\mathbb{L}=\mathbb{Q}$ of Theorem 2.7, we can
find an elementary end extension $\mathrm{I}_{\mathbb{Q}}$ of $\mathrm{I}$
such that SSy$_{\mathrm{I}}$($\mathrm{I}_{\mathbb{Q}})=\mathcal{A}$ and $\mathrm{I}_{%
\mathbb{Q}}\backslash \mathrm{I}$ contains a copy of $\mathbb{Q}^{'}:=\left\{ c_{q}:q\in
\mathbb{Q}\right\} $ of $\mathbb{Q}$, and there is a composition preserving
embedding\textit{\ }$j\mapsto \widehat{j}$ from the semi-group of initial
self-embeddings of $\mathbb{Q}$ into the\textit{\ }semi-group of initial
self-embeddings of $\mathrm{I}_{\mathbb{Q}}$ that satisfies conditions (a)
through (d) of Theorem 2.7. In particular $\mathbb{Q}^{'} $ is downward cofinal in $\mathrm{I}_{\mathbb{Q}}\backslash \mathrm{%
I.}$\medskip

\textbf{Stage 3:} An initial embedding $k:\mathcal{N}\rightarrow \mathrm{I}_{%
\mathbb{Q}}$ is constructed such that $k$ fixes each
element of $\mathrm{I}$. Note that Theorem 2.5 cannot be invoked for this
purpose since I$\Sigma _{1}$ need not hold in $\mathcal{N}$; instead, we
will take advantage of recursive saturation of $\mathcal{N}$, and the
properties of $\mathcal{N}$ listed in Remark 2.2. We will go through construction of $ k $ after describing stage 4 of the proof.\medskip

\textbf{Stage 4: }The desired self-embedding $j$ satisfying (3) of Theorem B
can be readily constructed as follows: Fix some $c_{q_{1}}<k(b_{0})$ in $%
\mathbb{Q}^{^{\prime }}$ and let $j_{1}$ be a fixed-point free initial
embedding of $\mathbb{Q}$ such that $j_{1}(\mathbb{Q})<q_{1}$. Then define $%
h:=k^{-1}\widehat{j_{1}}k$, and let $j$ be the restriction of $h$ to $%
\mathcal{M}$. First, note that by the way $j_{1}$ is chosen, $h$ is
well-defined and $h(\mathrm{N})<b_{0}$. Therefore, $j$ is an isomorphism
between $\mathcal{M}$ and a proper cut $\mathrm{J}$ of $\mathcal{M}$. Moreover, Since
$\mathrm{Fix}(\widehat{j_{1}})=\mathrm{I}$ (as arranged in Stage 2) and $k$
fixes each element of $\mathrm{I}$ (as arranged in Stage 3), by Remark 2.3
we may conclude that $\mathrm{Fix}(j)=\mathrm{I}$ and $j$ is an isomorphism
between $(\mathcal{M},\mathcal{A})$ and $(\mathrm{J},\mathcal{A}_{\mathrm{J}%
})$.\medskip

The above description of the four stages of the proof should make it clear
that the proof of condition (3) of Theorem B will be complete once we verify
that Stage 3 can be carried out, so we focus on the construction of an
initial embedding $k$ of $\mathcal{N}$ into $\mathrm{I}_{\mathbb{Q}}$ that
fixes each element of $\mathrm{I}$. To do so, we first note that since $(i)$
$\mathrm{I}\Sigma _{1}$ holds in both $\mathcal{M}$ and $\mathrm{I}_{\mathbb{%
Q}}$ , $(ii$) $\mathrm{SSy}_{\mathrm{I}}(\mathcal{M})=\mathrm{SSy}_{\mathrm{I%
}}(\mathrm{I}_{\mathbb{Q}})$, and $(iii)$ $\mathrm{Th}_{\Sigma _{1}}(%
\mathcal{M},i)_{i\in \mathrm{I}}=\mathrm{Th}_{\Sigma _{1}}(\mathcal{N}_{%
\mathbb{Q}},i)_{i\in \mathrm{I}}$ (because ${\mathrm{I}\prec _{\Sigma _{1}}%
\mathcal{M}}$, and ${\mathrm{I}}\prec \mathrm{I}_{\mathbb{Q}}$) by Theorem
2.5 there is a proper initial embedding $f:\mathcal{M}\rightarrow \mathrm{I}%
_{\mathbb{Q}}$ such that $f(i)=i$ for each $i\in I$ . In particular, $f(%
\mathrm{M})<e$ for some $e\in \mathrm{I}_{\mathbb{Q}}$. Moreover, since ${%
\mathrm{Th}_{\Pi _{1}}(\mathcal{M},x)_{x\leq b_{0}}=\mathrm{Th}_{\Pi _{1}}(%
\mathcal{N},x)_{x\leq b_{0}}}$ we have:

\begin{center}
$(\ast _{0}):\ \ \ \mathcal{N}\models \exists z\ \delta (z,w)\Rightarrow
\mathrm{I}_{\mathbb{Q}}\models \exists z<e\ \delta (z,f(w))$, for all $\Delta
_{0}$-formulas $\delta $ and all $w<b_{0}.$
\end{center}

\noindent Now choose $\{b_{n}:n\in \omega \}$ to be a decreasing sequence in $\mathrm{M}\setminus \mathrm{I}$ such that $b_{0}$ is the element
chosen in Stage 1, and $2^{b_{n+1}}<b_{n}$ for all $n\in \omega $. In order
to construct $k$, we recursively build finite sequences $\bar{u}%
:=(u_{0},...,u_{m})$ of elements of $\mathrm{N}$ and $\bar{v}%
:=(v_{0},...,v_{m})<e$ for each $m\in \omega $ such that:

\begin{itemize}
\item[$(i)$] $u_{0}=0=v_{0}$.

\item[$(ii)$] For every $c$ in $\mathrm{N}$ there is some $n\in \omega $
such that $c=u_{n}$.

\item[$(iii)$] For every $n\in \omega $, $v_{n}<b$, and if for some $c$ in $%
\mathrm{I}_{\mathbb{Q}} $ it holds that $c<v_{n}$, then there is some $m\in \omega
$ such that $c=v_{m}$.

\item[$(iv)$] For every $m\in \omega $ the following condition holds:

$(\ast _{m}):\ \ $ $\mathcal{N}\models \exists z\ \delta (z,w,\bar{u}%
)\Rightarrow \mathrm{I}_{\mathbb{Q}}\models \exists z<e\ \delta (z,f(w),\bar{%
v})$ for all $\Delta _{0}$-formulas $\delta $ and all $w<b_{m}$
\end{itemize}

Let $\{a_{n}:n\in \omega \}$ and $\{d_{n}:n\in \omega \}$ respectively be
enumerations of element of $N$ and $\underline{e}\subset \mathrm{I}_{\mathbb{%
Q}}$, and $\langle \delta _{r}:r\in \mathrm{M}\rangle $ be a canonical
enumeration of $\Delta _{0}$-formulas in $\mathcal{M}$. The first step of
induction holds thanks to $(\ast _{0})$ and the choice of $u_{0}$ and $v_{0}$
in statement $(i)$. Next, suppose $\bar{u}:=(u_{0},...,u_{m})\in \mathrm{N}$
and $\bar{v}:=(v_{0},...,v_{m})<e$ are constructed, for given $m\in \omega $%
. In order to build $u_{m+1}$ and $v_{m+1}$ we distinguish two cases, one
handling the `forth' step and the other handling the `back' step of the
back-and-forth construction: \medskip

\textbf{CASE 0 ($\mathbf{m+1=2k}$, for some $\mathbf{k}\in \omega $):} In
this case if $a_{k}$ is one of elements of $\bar{u}$, put $u_{m+1}=u_{m}$
and $v_{m+1}=v_{m}$. Otherwise, put $u_{m+1}=a_{k}$ and define:

\begin{center}
$\mathrm{A}:=\lbrace\langle r,w\rangle<b_{m+1}: \mathcal{N}\models \exists z
\ \mathrm{Sat}_{\Delta_{0}}(\delta_{r}(z,w,\bar{u},a_{k}))\rbrace .$
\end{center}

\noindent Note that in the above definition, $\mathcal{N}$ can be safely
replaced by some $\mathcal{N}_{n}$, where $n$ is large enough to contain the
parameters $\bar{u}$ and $a_{k}$ (thanks to properties $(i)$ and $(ii)$ in
Remark 2.2). On the other hand, by property $(iii)$ in Remark 2.2, there is
some $\mathrm{ED}_{n}\in \mathcal{A}$ such that:

\begin{center}
$\mathrm{A}=\left\{ \langle r,w\rangle <b_{m+1}:\ulcorner \exists z\ \mathrm{%
Sat}_{\Delta _{0}}(\delta _{r}(z,w,\bar{u},a_{k})\urcorner \in \mathrm{ED}%
_{n}\right\} .$
\end{center}

\noindent Since $(\mathcal{M},\mathcal{A})$ satisfies $\mathrm{I}\Sigma
_{1}^{0}$, the above characterization of $\mathrm{A}$ shows that $\mathrm{A}$
is coded in $\mathcal{N}$ by some $a<2^{b_{m+1}}.$ Therefore we have:\medskip

\noindent (1) $\mathcal{N}\models \forall \langle r,w\rangle
<b_{m+1}(\langle r,w\rangle Ea\rightarrow \exists z\ \mathrm{Sat}_{\Delta
_{0}}(\delta _{r}(z,w,\bar{u},a_{k}))).$\medskip

\noindent Recall that $\mathrm{B}\Sigma _{1}$ holds in $\mathcal{N}$, and $%
\mathrm{Sat}_{\Delta _{0}}$ has a $\Sigma _{1}$-description in $\mathcal{N},$ so (1) allows us to conclude:\medskip

\noindent (2) $\mathcal{N}\models \exists t\ \forall \langle r,w\rangle
<b_{m+1}(\langle r,w\rangle Ea\rightarrow \exists z<t\ \mathrm{Sat}_{\Delta
_{0}}(\delta _{r}(z,w,\bar{u},a_{k}))).$\medskip

\noindent By quantifying out $a_{k}$ in (2) and coupling it with $(\ast
_{m}) $, we obtain:\medskip

\noindent (3) $\mathrm{I}_{\mathbb{Q}}\models \exists x,t<e\ \forall \langle
r,w\rangle <f(b_{m+1})\ (\langle r,w\rangle Ef(a)\rightarrow \exists z<t\
\mathrm{Sat}_{\Delta _{0}}(\delta _{r}(z,w,\bar{v},x))).$\medskip

\noindent Clearly any element of $\mathrm{I}_{\mathbb{Q}}$ that witnesses $x$
in (3) can serve as a suitable candidate for $v_{m+1}$.\medskip

\textbf{CASE 1 ($\mathbf{m+1=2k+1}$, for some $\mathbf{k}\in \omega $):} In
this case if ${d_{k}\geq \mathrm{Max}\{\bar{v}\}}$ or if it is one of the
elements of $\bar{v}$, put $u_{m+1}=u_{m}$ and $v_{m+1}=v_{m}$. Otherwise,
put $v_{m+1}=b_{k}$ and define:

\begin{center}
$\mathrm{B}:=\{\langle r,w\rangle <f(b_{m+1}):\mathrm{I}_{\mathbb{Q}}\models
\forall z(\mathrm{Sat}_{\Delta _{0}}(\delta _{r}(z,w,\bar{v}%
,d_{k}))\rightarrow b<z)\}$.
\end{center}

\noindent Note that $\mathrm{B}$ is $\Sigma _{1}$-definable in $\mathrm{I}_{%
\mathbb{Q}}$, so there is some $b<2^{f(b_{m+1})}=f(2^{b_{m+1}})$ which codes
$\mathrm{B}$ in $\mathrm{I}_{\mathbb{Q}}$. Therefore $b=f(c)$ for some $%
c<2^{b_{m+1}}$. Define:

\begin{center}
$p(x):=\{\forall w<b_{m+1}(\langle \ulcorner \delta \urcorner ,w\rangle
Ec\rightarrow \forall z\ \lnot \delta (z,w,\bar{u},x)):\delta \ \mathrm{is\ a%
}\ \Delta _{0}\mathrm{-formula}\}$.
\end{center}

\noindent Since $\mathcal{N}$ is recursively saturated and $p(x)$ is
recursive, in order to find a suitable element in $\mathrm{N}$ which serves
as $u_{m+1}$ , it suffices to prove that $p(x)$ is finitely satisfiable. So
suppose $p(x)$ is not finitely satisfiable. It can be readily checked that $%
p(x)$ is closed under conjunction, so we can safely assume there is a $%
\Delta _{0}$-formula $\delta $ such that:\medskip

\noindent (4) $\mathcal{N}\models \forall x(\exists w<b_{m+1}(\langle
\ulcorner \delta \urcorner ,w\rangle E c \wedge \exists z \ \delta (z,w,\bar{%
u},x))).$\medskip

\noindent Clearly (4) implies:\medskip

\noindent (5) $\mathcal{N}\models \forall x<\mathrm{\mathrm{Max}}\lbrace\bar{%
u}\rbrace \ (\exists w<b_{m+1}(\langle \ulcorner \delta \urcorner ,w\rangle
E c \wedge \exists z \ \delta (z,w,\bar{u},x))).$\medskip

\noindent We can bound variable $z$ in (5) by using $\mathrm{B}\Sigma _{1}$
in $\mathcal{N}$, and next employ $(\ast _{m})$ to deduce:\medskip

\noindent (6) $\mathrm{I}_{\mathbb{Q}}\models {\exists t<e\forall x<\mathrm{%
\mathrm{Max}}\{\bar{v}\}(\exists w<f(b_{m+1})(\langle \ulcorner \delta
\urcorner ,w\rangle Ef(c)\wedge \exists z<t\ \delta (z,w,\bar{v},x)))}.$%
\medskip

\noindent By replacing $x$ in (6) with $d_{k}$, we obtain:\medskip

\noindent (7) $\mathrm{I}_{\mathbb{Q}}\models \exists t<e\ (\exists
w<f(b_{m+1})(\langle \ulcorner \delta \urcorner ,w\rangle Eb\wedge \exists
z<t\ \delta (z,w,\bar{v},d_{k}))).$\medskip

\noindent But (7) contradicts the assumption that  $b$ codes $\mathrm{%
B}$ in $\mathrm{I}_{\mathbb{Q}}$. So $p(x)$ is finitely satisfiable.\medskip
\hfill $\square $

\end{document}